\documentclass[12pt]{article}

\usepackage{graphicx}
\usepackage{bm}% bold math
\usepackage[mathlines]{lineno}% Enable numbering of text and display math
\usepackage{amsfonts} 

\newcommand{\beq}{\begin{equation}}
\newcommand{\eeq}{\end{equation}}

\newcommand{\ba}{\begin{array}}
\newcommand{\ea}{\end{array}}

\newcommand{\ben}{\begin{enumerate}}
\newcommand{\een}{\end{enumerate}}

\newtheorem{theorem}{Theorem}
\newcommand{\bth}{\begin{theorem}}

\newtheorem{definition}{definition}
\newcommand{\bdf}{\begin{definition}}
\newcommand{\edf}{\end{definition}}

\newtheorem{notation}{notation}
\newcommand{\bnt}{\begin{notation}}
\newcommand{\ent}{\end{notation}}

\newtheorem{prop}{Proposition}
\newcommand{\bprop}{\begin{prop}}
\newcommand{\eprop}{\end{prop}}

\newtheorem{lemma}{Lemma}
\newcommand{\blemma}{\begin{lemma}}
\newcommand{\elemma}{\end{lemma}}

\newtheorem{rmk}{Remark}
\newcommand{\brmk}{\begin{rmk}}
\newcommand{\ermk}{\end{rmk}}

\newtheorem{corr}{Corollary}
\newcommand{\bcorr}{\begin{corr}}
\newcommand{\ecorr}{\end{corr}}

\usepackage{xcolor}

\begin{document}
\title{On the Burgers dynamical system with an external force and its Koopman decomposition}
\author{Mikhael Balabane
\thanks{Laboratoire Analyse, G\'eom\'etrie et Applications, Universit\'e Paris 13}\thanks{Center for Advanced Mathematical Sciences, American University of Beirut} }
%- {Miguel Alfonso Mendez} \thanks{von Karman Institute for Fluid Dynamics, EA Dept, Sint Genesius Rode, Belgium} - {Sara Najem}
%\thanks{Physics Department, American University of Beirut, Beirut 1107 2020, Lebanon}}\date{oct 2020}
 \maketitle
 
 \begin{abstract}
 \noindent We prove that the Burgers flow with a steady external forcing has a unique steady state which is a sink. Although this flow cannot be linearized  through Cole-Hopf transforms, we prove that it has a convergent Koopman Modes decomposition. This gives an asymptotic formula for solutions of the Burgers equation with an external force. Time dependence and the coefficients of the decomposition are proved to be eigenvalues and eigenfunctionals of the Koopman operator. 

\noindent Convergence of the Koopman decomposition is proved for orbits close to the sink.

\noindent The analysis of Burgers dynamical system  relies on the properties of a nonlinear heat flow, that shows invariant sets  with complete orbits, and invariant sets where blow-up in finite time do ocurr. This behaviour helps understanding  some instabilities in numeric computing for  fluids.
 \end{abstract}

 \section{introduction}\label{intro}
Koopman flows were introduced by B.O.Koopman [1] to analyse hamiltonian finite dimensional flows. An important community in fluid dynamics makes use of the Koopman decomposition to analyse data originated from fluids computing. The Koopman decomposition for the field $u$ of velocities of a fluid is a representation of $u$ of the following form:
$$u(t,x)=\sum_\nu e^{-\lambda_\nu t}\varphi_\nu(u_0)a_\nu(x)$$
where $\nu$ belongs to some set of indices and $u_0$ is the initial condition. The $\lambda_\nu$'s are called the Koopman eigenvalues, the $a_\nu(x)$'s are called the Koopman modes.

%\noindent For a Hamiltonian flow $\Phi^t$ the Koopman flow $\mathcal{K}^t$ is a linear flow defined by Koopman [1] in order to 'linearise' finite dimensional Hamiltonian flows. It is a set of maps $(\mathcal{K}^t)_{t\in \mathbb{R}_+}$, where the index $t$ is time. Each map $\mathcal{K}^t$ maps observables ( i.e. continuous functionals $\phi$ of the state variable $u$) to $\mathbb{R}$, and  fulfills the following property:

\noindent In the important litterature on Koopman flows (N.J.Kutz - S.L.Brunton - B.W.Brunton - J.L.Proctor [3], I.Mezic [5], C.W.Rowley - S.T.Dawson [6], P.J.Schmid [7],..) proofs are only given for finite dimensional sytems of ODEs, but data analysis is performed for Navier-Stokes simulations. The only attempt for a PDE is, up to our knowledge,  J.Page and R.R.Kerswell's [4], where a procedure is given for the Burgers' equation with no external for\-cing. Their procedure could not be completed because they encountered a 'degeneracy' due to multiplicity of the eigenvalues of the Koopman operator. This degeneracy is overcome in M.Balabane-M.Mendez-S.Najem [2].

\noindent This paper is devoted to the  dynamical system defined by the Burgers equation on $[0,1]$, with an  external forcing $2\partial_xV$:
\begin{equation}\label{burgers}\partial_tu+u\partial_xu=\partial^2_{xx}u+2\partial_xV
\end{equation}
$$u(t,0)=u(t,1)=0\quad u(0,x)=u_0(x)
$$
 We prove that the dynamical system associated with this PDE has a unique sink, attracting all orbits. Notice that the function $V$ in (\ref{burgers}) is defined up to an additive constant, so one can assume $V>0$.
 
 \noindent Let $\Phi^t_\mathcal{B}(u_0)=u(t,.)$. A Koopman decomposition for $\Phi^t_\mathcal{B}$ is proved. We prove that the time dependences  in the decomposition are given by the eigenvalues of the Koopman operator $\mathcal{K}^t_\mathcal{B}$, and the coefficients are values at the Cauchy data of the eigenfunctionals of the Koopman operator, which is defined (Koopman [])  by:
 $$(\mathcal{K}^t_{\mathcal{B}} \phi))(u_0)=\phi(\Phi^t_{\mathcal{B}}(u_0))
$$
where $\phi$ is any continuous functional on $L^2([0,1])$, called observable in the DMD community.

\noindent The proofs are done through the Cole-Hopf transforms, although in the pre\-sence of the forcing this transform do not intertwine Burgers equation with a linear equation, but with a nonlinear and nonlocal heat equation, namely:
\begin{equation}\label{Nheat}\partial_t\tilde{v}=(\partial^2_{xx}-V) \tilde{v}+(\int_0^t V(x)\tilde{v}(t,x)dx\,)\tilde{v}
\end{equation}
$$\partial_x\tilde{v}(t,0)=\partial_x\tilde{v}(t,1)=0\quad \tilde{v}(0,x)=\tilde{v}_0(x)
$$
Properties of the flow associated with (\ref{Nheat}) are proven, and its Koopman decomposition is derived.

\noindent The behaviour of this flow gives hints for instabilities in fluid computations.
\medskip

\noindent To prove properties of this flow one relies on precise estimates for the fol\-lowing  linear  heat flow:
 \begin{equation}
 \partial_t {v}  =\partial_{xx}^2{v}  -V(x){v} \quad \partial_x {v}(t,0)=\partial_x {v}(.,1)=0\quad
{v}(0,x)= {v}_0(x)\in \mathcal{P}_+
 \label{eqH}
 \end{equation}
\noindent  In Section \ref{sCH} the Cole-Hopf transforms are given. In section \ref{slinks} links between the flows are proved. In section  \ref{sHfs}  the properties of the dynamical systems defined by   heat equations are proved: in  subsection  \ref{sslHf}  are overviewed the needed basic properties of the linear heat flow (\ref{eqH}); in subsection \ref{ssNf} are stated the properties of  the nonlinear heat flow defined by  (\ref{Nheat}).  Section  \ref{sBurgers}  is devoted to Burgers dynamical system  (\ref{burgers}) and its Koopman decomposition. Straightforward proofs are given in the core of the text, and lengthy proofs are  given in section  \ref{proofs}.

%-----------------------------------------------------------------------

\section{the Cole-Hopf transforms}\label{sCH}
The Cole-Hopf transforms are given by:
\begin{equation}\label{CH}
u:=C(v)=-2{\partial_xv\over v}\quad {\rm and} \quad v:=H(u)={e^{-{1\over 2}\int_0^xu(s)ds}\over \int_0^1e^{-{1\over 2}\int_0^xu(s)ds}dx}
\end{equation}

\noindent $H$ is defined on $L^2([0,1])$. $H(u)$ is  strictly positive on $[0,1]$, square integrable as well as its first weak derivative on $[0,1]$, and fulfills $\int_0^1H(u)(s)ds=1$. The set of  functions with these properties is denoted by $\mathcal{P}_1$. 

\noindent $C$ is defined for strictly positive functions on $[0,1]$, square integrable as well as their first weak derivative (hence continuous). The set of  functions fulfilling these properties is denoted by $\mathcal{P_+}$. 

\noindent One has, for all $ u\in L^2([0,1]), C(H(u))=u$.     If $v\in \mathcal{P}_1$  then $H(C(v))=v$. 

\noindent It follows that $H$ is a bijection from $L^2([0,1])$ onto $\mathcal{P}_1$, with inverse $C$. Both maps are differentiable.

\noindent A useful remark is that for any (spatially) constant $\delta$ one has 
\begin{equation}\label{inv}C(\delta v)= C(v)\end{equation}
Conversely if $v_1$ and $v_2$ belong to $  \mathcal{P}_+$ and fulfill $C(v_1)=C(v_2)$ then $g(v_1)v_2(x)=g(v_2) v_1(x)$ where $g(v)=\int_0^1v(s)ds$. Hence $v_1\over v_2$ do not depend on the space variable $x$.

%-------------------------------------------------------------------
\section{linking the flows}\label{slinks}
%--------------------------------------------------------
\subsection{Intertwining  ${\Phi^t_{\mathcal{B}}}$ and ${\Phi^t_{\mathcal{N}}}$ }\label{ssBN}
Let $\Phi^t_{\mathcal{B}}$ denote  the flow associated with Burgers equation (\ref{burgers}), meaning that $u(t,.)={\Phi^t_{\mathcal{B}}}(u_0)$ solve Burgers equation (\ref{burgers}) with Cauchy data $u_0\in L^2([0,1])$.  $\Phi^t_{\mathcal{B}}$ is a flow in $L^2([0,1])$.

\noindent  Let $\tilde{v}_0=H(u_0)$ and $\tilde{v}(t,.)=H(u(t,.))$, so that $u(t,.)=C(\tilde{v}(t,.))$. Then
$$\forall t\geq 0,\, \tilde{v}(t,.)\in \mathcal{P}_1\quad 
\tilde{v}(0,.)=\tilde{v}_0,\,\, \partial_x\tilde{v}(t,0)=\partial_x\tilde{v}(t,1)=0,\,\, \int_0^1\tilde{v}(t,\xi)d\xi=1
$$
In order to get a PDE for $\tilde{v}$ one notices that:
$$\partial_tu(t,.)=-2\partial_t{\partial_x\tilde{v}(t,.)\over \tilde{v(t,.)}}=-2\partial_x{\partial_t\tilde{v}(t,.)\over \tilde{v(t,.)}}
$$
so
$$\partial_tu+u\partial_xu=-2\partial_x[{\partial_t\tilde{v}\over\tilde{v}}-({\partial_x\tilde{v}\over\tilde{v}})^2]=
-2\partial_x[{\partial_t\tilde{v}\over\tilde{v}}-{\partial^2_{xx}\tilde{v}\over\tilde{v}}+\partial_x({\partial_x\tilde{v}\over \tilde{v}})]
$$
On the other hand
$$\partial_x(\partial_xu+2V)=-2\partial_x(\partial_x({\partial_x\tilde{v}\over \tilde{v}})-V)
$$
Equation (\ref{burgers}) is therefore equivalent to the existence of $n(t)$ such that:
$$\partial_t\tilde{v}=\partial^2_{xx}\tilde{v}-V\tilde{v}+n(t)\tilde{v}
$$

\noindent Integrating this equality gives:
$$n(t)=\int_0^1V(x)\tilde{v}(t,s)ds
$$
and proves that $\tilde{v}$ solves the nonlinear heat equation (\ref{Nheat}). 

\noindent The converse is true because the computation is performed through equivalences.  

\noindent This proves that if $\Phi^t_\mathcal{N}$ denotes the flow associated with equation (\ref{Nheat}) then:
\begin{equation}\label{iBN}
H\circ {\Phi^t_{\mathcal{B}}}={\Phi^t_{\mathcal{N}}}\circ H \quad {\rm and}\quad C\circ {\Phi^t_{\mathcal{N}}}={\Phi^t_{\mathcal{B}}}\circ C
\end{equation}
 In sections \ref{sssNpos} and \ref{sssNmean}  $\Phi^t_{\mathcal{N}}$ is proved to be a flow in $\mathcal{P}_1$.

%---------------------------------------------------------
\subsection{linking $\Phi^t_{\mathcal{N}}$ and $\Phi^t_{\mathcal{H}}$}\label{ssNC}
Let  $\Phi^t_\mathcal{H}$ denote the flow defined by equation (\ref{eqH}). It is a flow on $\mathcal{P}_+$, as proved in section \ref{sssHpos}. Let $v_0(x)=\tilde{v}_0(x)\in \mathcal{P}_1$. Let $v(t,.)=\Phi^t_{\mathcal{H}}(v_0)$ and $\tilde{v}(t,.)=\Phi^t_{\mathcal{N}}(\tilde{v}_0)$ for $t\in [0,T_*(\tilde{v_0})[$, the time span of $\Phi^t_{\mathcal{N}}(\tilde{v}_0)$.

\noindent Let $g(t)=\int_0^1v(t,\xi)d\xi$ so $g'(t)=-\int_0^1V(\xi)v(t,\xi)d\xi$. 

\noindent Let $f(t,x)=g(t)\tilde{v}(t,x)-v(t,x)$. 

\noindent Then:
$$\partial_tf-(\partial_{xx}^2-V)f-\tilde{v}(t,x)(\int_0^1V(\xi)f(t,\xi)d\xi)=0
$$
and
$$\partial_xf(t,0)=\partial_xf(t,1)=0
$$
This implies, for any $t\in [0,T]$ with $T<T_*(\tilde{v_0})$:
$${d\over dt}{1\over 2}\Vert f\Vert^2_{L^2}\leq \int_0^1V(\xi)f(t,\xi)d\xi\int_0^1\tilde{v}(t,\xi)f(t,\xi)d\xi
$$
$$\leq \Vert  V\Vert_{L^2} \sup_{[0,T]} \Vert \tilde{v} \Vert_{L^2} \Vert f \Vert^2_{L^2}
$$
But $f(0,x)=0$ so Gromwall's lemma gives $f=0$ and:
\begin{equation}\label{gvtilde}g\tilde{v}=v
\end{equation}
Because $\Phi^t_{\mathcal{H}}$ is a flow on ${\mathcal{P}_+}$ as recalled in section \ref{sslHf}, $g(t)>0$, so equation (\ref{gvtilde}) translates to:
\begin{equation}\label{lNC}
\Phi^t_{\mathcal{N}}={1\over g(t)}\Phi^t_{\mathcal{H}}
\end{equation}
Using (\ref{iBN}), and  invariance (\ref{inv}) of $C$ under a scaling,this gives for $t\in [0,T_*(\tilde{v_0})[$:
\begin{equation}\label{projtraj}
C\circ \Phi^t_{\mathcal{N}}=C\circ \Phi^t_{\mathcal{H}}= \Phi^t_{\mathcal{B}}\circ C
\end{equation}
Global existence of the flow of $\Phi^t_\mathcal{N}$ on $\mathcal{P}_1$ is proven in \ref{globN}, hence $T_*(\tilde{v_0})=\infty$
%-------------------------------------------------------------------
\section{on the linear  and nonlinear heat flows}\label{sHfs}
\subsection{on the linear  heat flow}\label{sslHf}
Assume $V$ is a regular, strictly positive function on $[0,1]$. Consider the  linear heat equation  (\ref{eqH}) on $[0,1]$ with $V$ for potential:
\begin{equation}\label{eqHbis}
 \partial_t {v}  =\partial_{xx}^2{v}  -V(x){v} \quad \partial_x {v}(t,0)=\partial_x {v}(t,1)=0\quad
{v}(0,x)= {v}_0(x)\in \mathcal{P}_+
\end{equation}
%-----------------------------------
\subsubsection{basic estimates}\label{sssHest}
Because $v_0$ is in $H^1$ one has the estimate:
$$\forall t\geq 0\quad \Vert v(t,.)\Vert_{H^1}\leq C \Vert v_0\Vert_{H^1}
$$
Regularising effect of the heat equation enables to extend this flow to Cauchy data in $L^2([0,1])$ and one has:
$$\forall t> 0\quad \Vert v(t,.)\Vert_{H^1}\leq {C_1\over \sqrt{t}} \Vert v_0\Vert_{L^2}
$$
Notice that the regularising effect extends to all derivatives, so the solution with Cauchy data in $L^2([0,1])$ lays in  all Sobolev spaces and one has the following estimates:
 $$\forall t> 0\quad \Vert v(t,.)\Vert_{H^s}\leq {C_s\over \sqrt{t^s}} \Vert v_0\Vert_{L^2}
$$
%---------------------------------------
\subsubsection{on positivity}\label{sssHpos}
 By Kato's lemma, if $v_0>0$ then $v(t,x)>0$ for all $t>0$, so the associated flow acts on $\mathcal{P}_+$. where it is denoted by  by $\Phi^t_{\mathcal{H}}$:
\begin{equation}
    \label{fH}
\forall t\geq 0,\forall {v}_0 \in \mathcal{P}_+, \quad \Phi^t_{\mathcal{H}}({v_0})={v}(t,.)>0
\end{equation}
We need below  an enhancement of this property, namely that:
\begin{equation}
    \label{minoration}
\forall t\geq 0,\forall x\in [0,1],\forall {v}_0 \in \mathcal{P}_+, \quad {v}(t,x)\geq e^{-t\sup V}\inf_{[0,1]}v_0(x)
\end{equation}It can be proved as follows: the heat equation with $V=0$ and Neumann boundary conditions is invariant by adding a constant to the state variable. Hence positivity preservation gives $e^{t\partial^2_{xx}}v_0\geq \inf_{[0,1]}v_0(x)$. The Trotter-Kato formula gives
$$\Phi^t_\mathcal{H}(v_0)(t,x)=\lim_{N\rightarrow \infty}(e^{{t\over N}\partial^2_{xx}}e^{-{t\over N}V(x)})^Nv_0(x)
$$
hence the lower bound. 
.%----------------------------------------------------
\subsubsection{on eigenvalues}\label{sssHeig}
The operator $A=-\partial_{xx}^2+V(x)$ is self-adjoint from $H^1(]0,1[)$ to its dual space, and has a bounded inverse. Hence by Rellich argument its spectrum is an increasing sequence of positive eigenvalues $(\mu_n)_{n\in \mathbb{N}}$ associated with a complete orthonormal set in $L^2([0,1])$ of eigenfunctions $(e_n(x))_{n\in N}$:
\begin{equation}
    \label{vpH}
\forall n \in {\mathbb{N}}, \quad \partial_{xx}^2e_n(x)-V(x)e_n(x)=-\mu_ne_n(x)\quad \partial_{x}e_n(0)=\partial_{x}e_n(1)=0
\end{equation}
\noindent Let:
\begin{equation}
    \label{cn}
\forall n \in {\mathbb{N}},  \forall {v}_0 \in H^1(]0,1[) \quad c_n(v_0)=\int_0^1e_n(x)v_0(x)dx
\end{equation}
then 
\begin{equation}
    \label{devH}v(t,x)=\sum_0^\infty e^{-\mu_nt}c_n(v_0)e_n(x)
    \end{equation}
Because $\mu_0$ is a variational extremum, $e_0(x)$ do not vanish in $]0,1[$. One can assume $e_0(x)>0$ on $]0,1[$. An orthogonality argument shows that it is the only eigenfunction with constant sign. One can even prove, comparing  the Neumann minimisation with the mixed Dirichlet-Neumann minimisation that  $e_0(x)>0$ on $[0,1]$. So 
\begin{equation}\label{e0>0}
m_0=\inf_{[0,1]}e_0(x)>0.
\end{equation}
\medskip

\noindent Note that  this implies: 
\begin{equation}
    \label{c0}
 \forall v_0\in \mathcal{P}_+\quad c_0(v_0)=m_0\int_0^1v_0(x)dx>0
\end{equation}
Note also that due to the variational caracterisation of $e_0$, multiplicity of $\mu_0$ is equal to one, therefore $\mu_0<\mu_1$.
%--------------------------------------------------------------
\subsubsection{on the spatial mean}\label{sssHmean}
The spatial mean of $v(t,x)$ solving (\ref{eqH}) is defined in \ref{ssNC} by: 
\begin{equation}\label{intva}
g(t)=\int_0^1v(t,\xi)d\xi
\end{equation}
Because $v_0\in \mathcal{P_+}$ positivity of $v$ implies
$$g'(t)=-\int_0^1V(x)v(t,x)dx\geq -(\sup{V})g(t)$$
so
\begin{equation} \label{intv}g(t)\geq e^{-(\sup{V})t}\int_0^1v_0(x)dx>0
\end{equation}
%------------------------------------------------------------------
\subsubsection{two asymptotic estimates}\label{sssHasymp}
\noindent The first estimate, proved in \ref{proof1}, is about the mean $g(t)$. It is needed  in \ref{sssKN} for the Koopman decomposition of solutions of the nonlinear heat equation (\ref{Nheat}). It states that for all $t\geq 0$, and all $v_0\in \mathcal{P}_+$, if  $v$ solves equation (\ref{eqH}), then:
\begin{equation}\label{estimation1}\vert e^{\mu_0 t}g(t)-c_0(v_0)\int_0^1e_0(\xi)d\xi \vert \leq e^{-(\mu_1-\mu_0)t}\Vert v_0-c_0(v_0)e_0(x)\Vert_{L^2} \Vert 1-c_0(1)e_0(x)\Vert_{L^2}
\end{equation}
\bigskip

\noindent  The second  estimate, proved in \ref{proof2}, is needed for the Koopman decomposition of the Bur\-gers flow (\ref{burgers})  given in section \ref{ssKBf}.  It states that:
\begin{equation}\label{estimation2}
 \forall t>0\quad \sup_{x\in [0,1]} \vert e^{\mu_0t}v(t,x)-c_0(v_0)e_0(x)\vert\leq \Vert v_0-c_0(v_0)e_0\Vert_{L^2}h_1(t)
\end{equation}
where $h_1(t)$ do not depend on $v_0$ (it is given explicitely in section \ref{ssKBf}). It is a decreasing function on $]0,\infty[$ with limit zero at infinity. It is given explicitely in section \ref{proof2}.
\bigskip

\noindent  Note that the two estimates are invariant by a multiplication of $v_0$ by a positive constant.

%--------------------------------------------------------------------
\subsubsection{on the Koopman flow}\label{sssHK}
\noindent Let  $\mathcal{F}_+$ be the set of continuous maps from $\mathcal{P}_+$ to ${\mathbb{R}}$.  The Koopman flow ${\mathcal{K}}^t_{\mathcal{H}}$  associated with $\Phi^t_{\mathcal{H}}$ acts on  $\mathcal{F}_+$  through:
\begin{equation}
    \label{fKH}
\forall \phi \in \mathcal{F}_+, \forall t\geq 0, \forall {v}_0 \in \mathcal{P}_+, \quad ({\mathcal{K}}^t_{\mathcal{H}}(\phi))({v}_0)=\phi(\Phi^t_{\mathcal{H}}({v}_0))
\end{equation}

\noindent For any $n\in \mathbb{N}$ the functional $c_n$ is an eigen-functional of the Koopman flow ${\mathcal{K}}^t_{\mathcal{H}}$, with eignevalue $e^{-\mu_n t}$, because for any $v_0 \in \mathcal{P}_+$:
\begin{equation}\label{vpKH1} ({\mathcal{K}}^t_{\mathcal{H}}(c_n))({v}_0)=c_n(\Phi^t_{\mathcal{H}}({v}_0))=c_n(\sum_{k\in N}e^{-\mu_k t}c_k(v_0)e_k(x))=e^{-\mu_n t}c_n(v_0)
\end{equation} 
Let $\nu=(q_0,q_1,..,q_m)$ with $q_i\in \mathbb{N}$. Let $\lambda_\nu=\sum_{i=0}^m\mu_{q_i}$.The multiplicative property of Koopman operators implies
:\begin{equation}
    \label{vpKH2}
  {\mathcal{K}}^t_{\mathcal{H}}\sigma_\nu=e^{-t\lambda_\nu}\sigma_\nu \quad  {\rm for}\quad\sigma_\nu:=\prod_{k=0}^mc_{q_k}
  \end{equation}
so $e^{-t\lambda_\nu}$ is an eigenvalue of ${\mathcal{K}}^t_{\mathcal{H}}$ for the eigen-functional $\sigma_\nu$.
\medskip

\noindent An interesting remark is that the Koopman decomposition of the flow $\Phi^t_{\mathcal{H}}$ is given by formula (\ref{devH}) where the only eigen-functionals showing  are those whose index is of length one. This is due to the linearity of the flow $\Phi^t_{\mathcal{H}}$ and the linearity of the functional $\delta_x$ that maps any function in $H^1$ to its value at location $x$ (note that on $H^1$ these are continuous functionals).

%--------------------------------------------------------------------------

\subsection{a nonlinear heat flow}\label{ssNf}
Cole-Hopf transforms intertwin Burgers equation with the following  nonlinear heat equation on $[0,1]$:

\begin{equation}\partial_t \tilde{v}  =\partial_{xx}^2 \tilde{v}  -V(x) \tilde{v} +(\int_0^1V(\xi) \tilde{v}(t,\xi)d\xi\,)\tilde{v} \label{eqN}
\end{equation}
$$
 \partial_x \tilde{v}(t,0)=\partial_x \tilde{v}(t,1)=0\quad
 \tilde{v}(0,x)= \tilde{v}_0(x)\in \mathcal{P}_+
 $$
 \noindent Local existence of the solution is granted in $L^2([0,1])$. Let $[0,T_*(\tilde{v}_0)[$ be the maximal time  span.
 %---------------------------------------------------------------------------
 \subsubsection{on positivity}\label{sssNpos}
The flow $\Phi^t_\mathcal{N}$ preserves po\-sitivity for $t\in [0,T_*(\tilde{v_0})[$. This can be checked using Duhamel formula and property (\ref{fH}) as follows, translating formula (\ref{eqN}) to:
$$\tilde{v}(t,x)=e^{t(\partial^2_{xx}-V)}\tilde{v}_0+\int_0^te^{(t-\tau)(\partial^2_{xx}-V)}[(\int_0^1V(\xi) \tilde{v}(\tau,\xi)d\xi\,)\tilde{v}(\tau, x)]d\tau
$$
As long as $\tilde{v}$ is positive, one has $(\int_0^1V(\xi) \tilde{v}(t,\xi)d\xi\,)\tilde{v}(t, x)\geq 0$, hence  applying section \ref{sssHpos},
$e^{(t-\tau)(\partial^2_{xx}-V)}[(\int_0^1V(\xi) \tilde{v}(\tau,\xi)d\xi)\tilde{v}(\tau, x)]\geq 0$, hence  
\begin{equation}\label{duhamel} 
\tilde{v}(t,x)\geq e^{t(\partial^2_{xx}-V)}\tilde{v}_0\geq e^{-t\sup V}\inf_{[0,1]}\tilde{v}_0(x)
\end{equation}
therefore  $\tilde{v}$ is trapped above  $e^{-t\sup V}\inf_{[0,1]}\tilde{v}_0(x)$ as long as it is positive, so:
\begin{equation}\label{tildev>0}
\forall t\in [0,T_*[, \forall x \in [0,1] \quad \tilde{v}(t,x)>0
 \end{equation}

\noindent Duhamel formula (\ref{duhamel}) shows moreover that the solution belongs to $H^1([0,1])$  for all $t\in [0,T_*(\tilde{v_0})[$.

\noindent Positivity preservation and regularity  imply 
that $\mathcal{P}_+$ is invariant by the flow  defined by equation (\ref{eqN}). 

%-------------------------------------------------------------------------

\subsubsection{on the spatial mean}\label{sssNmean}
\noindent Let $\tilde{v}_0 \in \mathcal{P}_+$. Let $\tilde{g}_0=\int_0^1\tilde{v}_0(\xi)d\xi$. Let $\tilde{g}(t)=\int_0^1\tilde{v}(t,\xi)d\xi $ where $ \tilde{v}(t,x)$  solves  (\ref{Nheat}) on its maximal time span $t\in [0,T_*(\tilde{v}_0)[$. One has
\begin{equation}\label{eqmean}{d\over dt}(\tilde{g}(t)-1)=(\int_0^1V(\xi)\tilde{v}(t,\xi)d\xi)\,\,(\tilde{g}(t)-1)
 \end{equation}
This implies for $t\in [0,T_*(\tilde{v_0})[$:

$$ \tilde{g}(t) > 1 \quad {\rm if} \quad \tilde{g}_0> 1  $$
\begin{equation}\label{intgtilde}  \tilde{g}(t) = 1\quad {\rm if} \quad \tilde{g}_0= 1 
 \end{equation}
$$\tilde{g}_0e^{-t\sup V}\leq \tilde{g}(t) < 1\quad  {\rm if} \quad 0\leq  \tilde{g}_0< 1  $$

\noindent Note that formula  (\ref{intgtilde}) and positivity of $\tilde{v}$ proven in \ref{sssNpos}   show through formula  (\ref{eqmean}) that $\tilde{g}$ is decreasing if $\tilde{g}_0<1$,  and increasing if $\tilde{g}_0>1$. 

\noindent Formula (\ref{tildev>0}) shows that positivity is preserved by the flow defined by equation (\ref{Nheat}). This added to formula (\ref{intgtilde})  shows that $\mathcal{P}_1$ is invariant by this flow. We denote  it by $\Phi^t_{\mathcal{N}}$ as acting in $\mathcal{P}_1$:
\begin{equation}
    \label{fN}
\forall t\geq 0,\forall \tilde{v}_0 \in \mathcal{P}_1, \quad \Phi^t_{\mathcal{N}}(\tilde{v_0})=\tilde{v}(t,.)\in \mathcal{P}_1
\end{equation}
\medskip
%---------------------------------------------------------------
\subsubsection{on global existence}\label{globN}
Let $0\leq \tilde{g_0}\leq 1$. For $t\in [0,T_*(\tilde{v_0})[$ integration by parts gives:
$${d\over dt}{1\over 2}\Vert \tilde{v}\Vert^2_{L^2}\leq (\int_0^1V(\xi)\tilde{v}(t,\xi)d\xi)\Vert \tilde{v}\Vert_{L^2}^2\
$$
Positivity proven in \ref{sssNpos} gives:
$${d\over dt}{1\over 2}\Vert \tilde{v}\Vert^2_{L^2}\leq (\sup_{[0,1]}{V})\tilde{g}(t)\Vert \tilde{v}\Vert_{L^2}^2\
$$
Formula (\ref{intgtilde}) implies:
$${d\over dt}{1\over 2}\Vert \tilde{v}\Vert^2_{L^2}\leq (\sup_{[0,1]}{V})\Vert \tilde{v}\Vert_{L^2}^2\
$$
This implies global existence of $\tilde{v}$ in $L^2([0,1])$ if $0\leq \tilde{g}_0\leq 1$, and the estimate:
\begin{equation}\label{exglob}
 \Vert \tilde{v}\Vert_{L^2}\leq e^{t\sup V}\Vert \tilde{v}_0\Vert_{L^2})
\end{equation}
Duhamel formula in \ref{sssNpos} grants global existence in $H^1([0,1])$. 
%-------------------------------------------------------------
\subsubsection{the case $\tilde{g}_0=1$:  global existence and regularity}\label{gtilde0=1}
If $\tilde{g}_0=1$ an important link between equation (\ref{Nheat}) and equation (\ref{eqH}) is the following: let $v$ solve (\ref{eqH}) with initial condition in $\mathcal{P}_1$. Because of formula (\ref{intv})  the following function $\tilde{w}={v\over g(t)}$ is well-defined, fulfills $\tilde{w}_0=\tilde{w}(0,x) \in \mathcal{P}_1$, and:
$$\partial_t\tilde{w}={\partial_t v\over g(t)}-{\partial_t g\over g^2}v={1\over g(t)}(\partial_{xx}^2-V)v+{(\int_0^1Vv)\over g^2}v=(\partial_{xx}^2-V)\tilde{w}+(\int_0^1V\tilde{w})\tilde{w}
$$
so $\tilde{w}$ solves equation (\ref{Nheat}). This implies
\begin{equation} \label{vetvtilde}\tilde{v}_0\in \mathcal{P}_1\Longrightarrow  \tilde{v}={v\over g}\quad {\rm with}\quad v_0=\tilde{v}_0
\end{equation}
Note that this formula gives another proof of global existence of the solution of (\ref{Nheat}), because $v(t,.)$ exists globally for all $t\geq 0$, and $g(t)$ is strictly positive for all $t\geq 0$, by formula (\ref{intv}).

\noindent This formula has also the following  important consequence: for $\tilde{v}_0$ in $\mathcal{P}_1$ the solution $\tilde{v}$ is infinitely differentiable, due to the regularity of the solution of the linear heat equation as stated in \ref{sssHest}
%-------------------------------------------------------------
\subsubsection{on blow-up}\label{blowup}

\noindent If $\tilde{g}_0> 1$ the solution of equation (\ref{Nheat}) blows-up in finite time.  This  follows the blow-up of $\tilde{g}$, which is a consequence of the following estimate:
$${d\over dt}\tilde{g}(t)=(\int_0^1V(\xi)\tilde{v}(t,\xi)d\xi)\,\,(\tilde{g}(t)-1)\geq (\inf_{[0,1]}V)\,\,(\tilde{g}(t)-1)\tilde{g}(t)
$$
\medskip

\noindent A complete description of the blow-up of $\tilde{v}$ can be given through the change of state variable:
 $$\tilde{z}={\tilde{v}\over \tilde{g}}$$
 where $\tilde{v}$ solves (\ref{Nheat}) with $\tilde{v}(0,x )\in \mathcal{P}_+$. $\tilde{z}$  fulfills:
$$\partial_t\tilde{z}={\partial_t \tilde{v}\over \tilde{g}}-{\partial_t \tilde{g}\over \tilde{g}^2}\tilde{v}=(\partial_{xx}^2-V)\tilde{z}+{(\int_0^1V\tilde{z})}\tilde{z}
$$
Therefor $\tilde{z}$ solves equation (\ref{Nheat}) and $\tilde{z}(0,x)\in \mathcal{P}_1$, hence $\tilde{z}$ exists for all $t\geq 0$ by \ref{globN}. Equality $\tilde{g}\tilde{z}=\tilde{v}$  implies that $\tilde{v}$ blows-up at all spatial locations at the same blow-up time:  the blow-up time of $\tilde{g}$.
\medskip

\noindent The behaviour of solutions of (\ref{Nheat}) on $\mathcal{P}_+$ can help understand some numrical instabilities of fluids computation, computing Burgers solution being, by a change of the state variable, computing solutions of (\ref{Nheat}) on $\mathcal{P}_1$.  
%----------------------------------------------------
\subsubsection{on steady states}\label{steadyN}

If $\tilde{s}(x)$ is a steady state solution of the nonlinear heat equation (\ref{Nheat})  in $\mathcal{P}_+$, it solves:
\begin{equation} \label{quasieq}\partial_{xx}^2 \tilde{s}  -V(x) \tilde{s} +(\int_0^1V(\xi) \tilde{s}(\xi)d\xi\,)\tilde{s}=0
\end{equation} 
with Neumann boundary conditions. Equation (\ref{quasieq}) states that $\tilde{s}$ is an eigenfunction of $\partial_{xx}^2   -V(x)$ that belongs to $\mathcal{P}_+$. Hence $\tilde{s}(x)=Ce_0(x)$ for some $C\in\mathbb{R}_+$ as proved in \ref{sssHeig}. Computing $C$ from (\ref{quasieq}) shows uniqueness of the steady state in $\mathcal{P}_+$, namely:
\begin{equation}\label{equasistate} \tilde{f}_0(x)={\mu_0\over \int_0^1V(\xi)e_0(\xi)d\xi}\,e_0(x)
\end{equation}
Integrating equation (\ref{vpH})  gives $\int_0^1V(\xi)e_0(\xi)d\xi=\mu_0\,\int_0^1e_0(\xi)d\xi$, so: 
\begin{equation} \label{quasistate}\tilde{f}_0(x)={e_0(x)\over \int_0^1e_0(\xi)d\xi}
\end{equation}
This proves that $\tilde{f}_0\in \mathcal{P}_1$, so the flow on  $\mathcal{P}_1$ has a unique steady state.
It is a sink: this can be checked using the formula (\ref{lNC}) and asymptotics (\ref{estimation1}) and (\ref{estimation2}) with $t\rightarrow \infty$ as follows:
$$\lim_{t\rightarrow \infty}\tilde{v}(t,x)=\lim_{t\rightarrow \infty}{e^{\mu_0t}v(t,x)\over e^{\mu_0t}g(t)}={c_0(v_0)e_0(x)\over c_0(v_0)\int_0^1e_0(\xi)d\xi}=\tilde{f}_0(x)
$$
%-------------------------------
\subsubsection{on spectral elements of the Koopman operator}\label{sssSKN}
\noindent In this section are computed spectral elements of the Koopman flow $\mathcal{K}^t_\mathcal{N}$ associated with the flowf $\Phi^t_\mathcal{N}$. These elements are the building blocks of the Koopman decomposition of $\Phi^t_\mathcal{N}$ performed in the next section.  
\medskip

\noindent To give a precise definition of the Koopman operator $\mathcal{K}_\mathcal{N}^t$, one needs the set $\mathcal{F}_1$ of observables on $\mathcal{P}_1$: it is the set of  maps from $\mathcal{P}_1$ to $\mathbb{R}$ which are continuous for the $H^1$ norm given by:
$$\Vert \tilde{v}_0\Vert_{H^1}^2=\Vert \tilde{v}_0\Vert_{L^2}^2+\Vert \partial_x \tilde{v}_0\Vert_{L^2}^2
$$
The Koopman operator $\mathcal{K}_\mathcal{N}^t$ maps $\phi \in \mathcal{F}_1$ to $\mathcal{K}_\mathcal{N}^t(\phi) \in \mathcal{F}_1$ by the formula:  
\begin{equation}\label{KON}\forall \tilde{v}_0\in \mathcal{P}_1\quad (\mathcal{K}_\mathcal{N}^t(\phi))(\tilde{v}_0)=\phi(\Phi^t_\mathcal{N}(\tilde{v}_0))
\end{equation}
This definition is consistent because $\mathcal{P}_1$ is invariant under $\Phi^t_\mathcal{N}$ due to (\ref{tildev>0}), (\ref{intgtilde}) and (\ref{exglob}).

\noindent To compute spectral elements of $\mathcal{K}_\mathcal{N}^t$ one notices  that  all functionals $c_n$ given by (\ref{cn}) belong to $\mathcal{F}_1$. Moreover  strict positivity of  $c_0$ on $\mathcal{P}_1$ proved in (\ref{c0}) shows that  $1\over c_0$ is an observable on $\mathcal{P}_1$. Therefore
\begin{equation}\label{cnsc0}\forall q\in \mathbb{N}, \quad{c_q\over c_0} \in \mathcal{F}_1
\end{equation}
Formulas (\ref{lNC}) and (\ref{vpKH1}) give, for $q\geq0$:
$$\forall \tilde{v}_0 \in \mathcal{P}_1\quad \mathcal{K}^t_\mathcal{N}({c_q})(\tilde{v}_0)=c_q(\Phi^t_\mathcal{N}(\tilde{v}_0))={1\over g(t)}c_q(\Phi^t_\mathcal{H}\tilde{v}_0))={1\over g(t)}e^{-\mu_q t}c_q(\tilde{v}_0)
$$
Therefore, the multiplicative property of $\mathcal{K}^t_\mathcal{N}$ gives:
\begin{equation}\label{vpKN1}\forall q\geq 0,\quad  \mathcal{K}^t_\mathcal{N}({c_q\over c_0})=e^{-(\mu_q-\mu_0) t}{c_q\over c_0}
\end{equation}
Applying once more the multiplicative property gives eigenvalues and eigen-observables of $\mathcal{K}^t_\mathcal{N}$, namely that $e^{-\lambda_\nu t }$ is an eigenvalue of  $\mathcal{K}^t_\mathcal{N}$, with eigenfunctional $\psi_\nu$, where:
\begin{equation} \label{vpKN2}
 \nu=(q_0,q_1,...,q_m)\in  \mathbb{N}\times( \mathbb{N}\backslash \{ 0 \})^m \quad {\rm for}\quad m\in  \mathbb{N}
\end{equation}
\begin{equation} \label{vpKN3}\lambda_\nu=\sum_{i=0}^k(\mu_{q_i}-\mu_0)\quad{\rm and}\quad
\psi_\nu=\prod_{i=0}^m{c_{q_i}\over c_{0}}
\end{equation}
%--------------------------------------------------
\subsubsection{on Koopman decomposition}\label{sssKN}
Let $\tilde{v}_0\in \mathcal{P}_1$. In order to motivate the Koopman decomposition for $\Phi^t_\mathcal{N}(\tilde{v}_0)$, it is worth  considering  formulas (\ref{vetvtilde}) and (\ref{devH}),  written using (\ref{e0>0}) as:
\begin{equation}\label{petitdevN}
\tilde{v}(t,x)=\Phi^t_\mathcal{N}(\tilde{v}_0)={1\over \int_0^1e_0(\xi)d\xi}\sum_{n=0}^\infty {e^{-(\mu_n-\mu_0)t}\over 1+ \tilde{k}(t)}{c_n(\tilde{v}_0)\over c_0(\tilde{v}_0)}e_n(x)
\end{equation}
where
$$\tilde{k}(t)=\sum_{q=1}^\infty {e^{-(\mu_q-\mu_0)t}}{c_q(\tilde{v}_0)\over c_0(\tilde{v}_0)}\,\,p_q\quad {\rm with}\quad p_q={ \int_0^1e_q(\xi)d\xi\over \int_0^1e_0(\xi)d\xi}
$$
Formula (\ref{petitdevN}) is not a Koopman decomposition for $\Phi^t_\mathcal{N}(\tilde{v}_0)$ because the time evolution is not given by exponentials. But all quatities involved are spectral elements of the Koopman operator, computed in the above section.
\medskip

\noindent In order to derive a Koopman decomposition for $\Phi^t_\mathcal{N}(\tilde{v}_0)$ one must expand $(1+\tilde{k})^{-1}=\sum_{m=0}^\infty(-1)^m\tilde{k}^m$. For this end one must have $\vert \tilde{k}\vert<1$.  Formula (\ref{estimation1}) gives:
$$\vert \tilde{k}(t)\vert\leq {e^{-(\mu_1-\mu_0)t}\over c_0(\tilde{v}_0)\int_0^1e_0(\xi)d\xi}\Vert \tilde{v}_0-c_0(\tilde{v}_0)e_0\Vert_{L^2} \Vert  1-c_0(1)e_0\Vert_{L^2}
$$
Let:
$$\tau_\mathcal{N}(\tilde{v}_0)={1\over \mu_1-\mu_0}\log({\Vert \tilde{v}_0-c_0(\tilde{v}_0)e_0\Vert_{L^2} \Vert  1-c_0(1)e_0\Vert_{L^2}
\over c_0(\tilde{v}_0)\int_0^1e_0(\xi)d\xi})
$$
Section \ref{proof3} gives the procedure by which, for $t\geq 0$ fulfilling $t> \tau_\mathcal{N}(\tilde{v}_0)$, formula (\ref{petitdevN}) can be written as:
\begin{equation}\label{KDN}\tilde{v}(t,x)=\sum_{q_0=0}^\infty  \sum_{m=0}^\infty\sum_{(q_1,..,q_m)\in( \mathbb{N}\backslash \{ 0 \})^m}e^{-\lambda_\nu t}\psi_\nu(\tilde{v}_0)b_\nu(x)
\end{equation}
with
$$\nu=(q_0,..,q_m)\in \mathbb{N}\times( \mathbb{N}\backslash \{ 0 \})^m\quad b_\nu(x)=(-1)^m{\prod_1^mp_{q_i}\over \int_0^1e_0(\xi)d\xi}{e}_{q_0}(x)\quad {\rm for}\quad m\neq 0
$$
and, for $m=0$, the last sum in (\ref{KDN}) reduced to one element, namely:  $$e^{-(\mu_{q_0}-\mu_0)t}{c_{q_0}(\tilde{v}_0)\over c_0(\tilde{v}_0)}{e_{q_0}(x)\over \int_0^1e_0(\xi)d\xi}$$

\noindent This is a Koopman decomposition for $\Phi^t_\mathcal{N}(\tilde{v}_0)$. The Koopman modes are the functions $ b_{\nu}(x)$. For each mode
the time behaviour is given by  the eigenvalue $e^{-\lambda_\nu t}$ of the Koopman operator. The coefficient is given by  the value of the associated eigen-functional $\psi_\nu$ at the Cauchy data. 

\noindent It is interesting to notice that each Koopman mode corresponds to (infinitely) many eigenvalues.

\noindent The only nondecreasing component in this decomposition is for $(q_0,m)=(0,0)$. It corresponds, as expected, to the sink $\tilde{f}_0$.

\noindent It is important to point out  that the Koopman decomposition is valid for all $t\geq 0$ (i.e. $ \tau_\mathcal{N}(\tilde{v}_0)<0$)  if the Cauchy data $\tilde{v}_0$ belongs to a neighbourhood of the sink $\tilde{f}_0$, namely:  
\begin{equation}\label{vicinN}\Omega_\mathcal{N}=\{\tilde{v}_0\in \mathcal{P}_1\,;\, {\Vert \tilde{v}_0-c_0(\tilde{v}_0)e_0\Vert_{L^2} \Vert 1-c_0(1)e_0\Vert_{L^2}
\over c_0(\tilde{v}_0)\int_0^1e_0(\xi)d\xi}<1\}
\end{equation}
Notice that although the spectral elements, building bricks of the Koopman decomposition, are defined on all  $\mathcal{P}_1$, the decomposition exists only locally on $\mathcal{P}_1$.
%-----------------------------------------------
\subsubsection{on absolute convergence}\label{sssAKN}
The series in formula (\ref{KDN}) is convergent under the assumption $\tilde{v}_0\in \Omega_\mathcal{N}$, but it is not absolutely convergent. This means that the summation order  is mandatory. In order to have  convergence  in the $\sup_x$ norm,  enabling through  Fubini theorem a summation in arbitrary order, and analytic use of the formula, through integration for instance, one needs a more restrictive assumption. To state it let:
\begin{equation}\label{S(t)}C_V=max(1,\sup_{x\in[0,1]}V(x))\quad{\rm and}\quad  h_2(t)=\sqrt{\sum_1^\infty {e^{-2(\mu_q-\mu_0)t}\over \mu_q^3}}
\end {equation}
The function $h_2(t)$ is defined for all $t\geq 0$, due to formula (\ref{mu3}) below. It is decreasing and has limit zero at infinity. 
For any $\tilde{v}_0 \in \mathcal{P}_1$ let:
\begin{equation}\label{tautilde}\tilde{\tau}_\mathcal{N}(\tilde{v}_0)=\inf_{t}\{t>0\,; \, h_2(t)<{c_0(\tilde{v}_0)\int_0^1e_0(\xi)d\xi \over C_V\Vert V\Vert_{L^2}\Vert \tilde{v}_0-c_0(\tilde{v}_0)e_0\Vert_{H^1}}\,\}
\end {equation}
If $t> \max(\tilde{\tau}_\mathcal{N}(\tilde{v}_0),\tau_\mathcal{N}(\tilde{v}_0))$ then 
the Koopman decomposition of $\tilde{v}=\Phi^t_\mathcal{N}\tilde{v}_0$ is given by the following series, which is absolutely convergent in the $\sup_x$ norm:
\begin{equation}\label{KN}\tilde{v}(t,x)=\sum_{\nu \in \mathcal{B}}e^{-\lambda_\nu t}\psi_\nu(\tilde{v}_0)b_\nu(x)
\end{equation}
where:
$$ \mathcal{B}=\cup_{m=0}^\infty \mathcal{B}_m$$
$$\forall m\in \mathbb{N}; \quad \mathcal{B}_m=\{\nu=(q_0,...,q_m); \,\,q_0\in {\mathbb{N}},\,\,{\rm and}\,\,q_i\in \mathbb{N}\backslash \{ 0\}\,\,{\rm for}\,\,i=1,..,m\}$$
The proof is given in section \ref{proof4}
\medskip

\noindent A consequence of formula (\ref{tautilde}) is that there exists a neighbourhood $\tilde{\Omega}_\mathcal{N}$ of the sink  $\tilde{f}_0$ such that for any $\tilde{v}_0 \in\tilde{\Omega}_\mathcal{N}$  absolute convergence occurs for all $t>0$. $\tilde{\Omega}_\mathcal{N}$  is  the set of functions $ \tilde{v}_0 \in \mathcal{P}_1$ such that:
\begin{equation}\label{Omegatilde}
\Vert \tilde{v}_0-c_0(\tilde{v}_0)e_0\Vert_{H^1}<      {c_0(\tilde{v}_0)\int_0^1e_0(\xi)d\xi \over {\max (C_V\Vert V\Vert_{L^2}h_2(0)}  , \Vert 1-c_0(1)e_0(x)\Vert_{L^2})}     
\end{equation}

%---------------------------------------
\section{on the Burgers flow}\label{sBurgers}
%----------------------------------------------------------------------
\subsection{on global existence and regularity}\label{ssregB}
The intertwining of the Burgers flow $\Phi^t_{\mathcal{B}}$ with the nonlinear heat flow $\Phi^t_{\mathcal{N}}$ through the Cole-Hopf transforms, as given by formula ({\ref{iBN}}), is valid for all initial data $u_0\in L^2([0,1])$ in the Burgers equation (\ref{burgers}), because $H$   maps bijectively $L^2([0,1])$ onto $\mathcal{P}_1$.

\noindent Therefore
the global existence of the nonlinear heat flow $\Phi^t_{\mathcal{N}}$ proved in \ref{gtilde0=1} implies  global existence for  $t\geq 0$ of the Burgers orbits $\Phi^t_{\mathcal{B}}(u_0)$, for all Cauchy data $u_0\in L^2([0,1])$.

\noindent It is worth noticing that the Cole-Hopf transform $C$ maps $(s+1)$-differentiable and strictly positive functions to $s$-differentiable functions. Therefore regula\-rity of the flow $\Phi^t_{\mathcal{N}}$ stated in \ref{gtilde0=1} implies that for any given Cauchy data $u_0\in L^2([0,1])$ the orbit $u(t,x)=\Phi^t_{\mathcal{N}}(u_0)$ lays in the set of infinitely differentiable functions in the space variable, for  $t>0$. Thus the productt $u\partial_xu$ in equation  (\ref{burgers}) is meaningful.
%----------------------------------------------------------------------
\subsection{on steady states for the Burgers flow}\label{sssteadyB}
The intertwining asserted by formula ({\ref{iBN}}), implies that a steady state for the Burgers flow $\Phi^t_{\mathcal{B}}$ is  image through $C$ of a steady state of the nonlinear heat flow $\Phi^t_{\mathcal{N}}$ on $\mathcal{P}_1$. Therefore section {\ref{steadyN}} shows that the only steady state for $\Phi^t_{\mathcal{B}}$ is the sink:
$$s_0=C(\tilde{f}_0)$$
Formula (\ref{inv}) gives
\begin{equation}\label{steadyB}
s_0(x)=C\big( {e_0(x)\over \int_0^1e_0(\xi)d\xi}\big)=-2{\partial_xe_0(x)\over e_0(x)}
\end{equation}
This can also be checked by a direct proof, that is given in section \ref{proofsink}.
%----------------------------------------------------------------------
\subsection{the linear Koopman flow for Burgers equation}\label{ssKBf}
Let $u_0\in L^2([0,1])$. Formula (\ref{projtraj}) gives the  Burgers orbit of $u_0$  as the image by $C$ of the linear heat orbit $\Phi^t_\mathcal{H}(v_0)$ where  $v_0=H(u_0)$ because:
$$\Phi^t_\mathcal{B}(u_0)=\Phi^t_\mathcal{B}\circ C(v_0)=C\circ \Phi^t_\mathcal{H}(v_0)=-2\,\,{\partial_x\Phi^t_\mathcal{H}(v_0)\over \Phi^t_\mathcal{H}(v_0)}
$$
Let $v(t,.)=\Phi^t_\mathcal{H}(v_0)$. The above formula gives:
$$u(t,.)=\Phi^t_\mathcal{B}(u_0)=-2{\partial_xv\over v}=-2{e^{\mu_0 t}\partial_xv\over {c_0(v_0)e_0+(e^{\mu_0 t}v-c_0(v_0)e_0)}}
$$
Estimate (\ref{estimation2}) gives:
$$
 \forall t>0\quad \sup_{x\in [0,1]} \vert e^{\mu_0t}v(t,x)-c_0(v_0)e_0(x)\vert\leq \Vert v_0-c_0(v_0)e_0\Vert_{L^2}h_1(t)
$$
so for 
\begin{equation}\label{KNass}{\Vert v_0-c_0(v_0)e_0\Vert_{L^2}\over c_0(v_0)m_0}h_1(t)\,<1
\end{equation}
one gets
$$u(t,.)=-2\,{e^{\mu_0 t}\partial_xv\over {c_0(v_0)e_0}}\,{1\over 1+ \tilde{k}_\mathcal{B}}=-2\,{e^{\mu_0 t}\partial_xv\over {c_0(v_0)e_0}}\sum_{m=0}^\infty (-1)^m\tilde{k}_\mathcal{B}^m
$$
where
$$
\tilde{k}_\mathcal{B}={e^{\mu_0 t}v-c_0(v_0)e_0\over c_0(v_0)e_0}
$$
Writing $v$ as the series given by formula (\ref{devH}), it is proved in \ref{proof4bis}, using the same procedure as for the Koopman flow of the nonlinear heat equation in \ref{sssKN}, that:
\begin{equation}\label{finalKF}u(t,x)=\sum_{q_0=0}^\infty \sum_{m=0}^\infty \sum_{(q_1,..,q_m)\in (\mathbb{N}\backslash\{ 0\})^m}^\infty e^{-\lambda_\nu t}\varphi_\nu(u_0)\,a_\nu(x)
\end{equation}
where
$$\forall m\in \mathbb{N},\, \forall \nu=(q_0,..,q_m) \in \mathcal{B}_m,\quad \lambda_\nu=\sum_0^m(\mu_{q_i}-\mu_0)$$
\begin{equation}\label{finuanu} \varphi_\nu(u_0)=\prod_0^m{ c_{q_i}(H(u_0))\over c_0(H(u_0))}\quad \quad a_\nu(x)=2(-1)^{m+1}{\partial_xe_{q_0}(x)\over e_0(x)}\prod_1^m{e_{q_i}(x)\over e_0(x)}
\end{equation}

\noindent Formula (\ref{finalKF}) is a Koopman decomposition for the Burgers flow. It needs  assumption (\ref{KNass}) be fulfilled, that can be rewritten in order to highlight the link with neighbourhoods of the sink $s_0$. let:
$$\forall u_0 \in L^2([0,1])\quad \tau_\mathcal{B}(u_0)=h_1^{-1}({m_0^2\over \Vert H(u_0)-\tilde{f}_0\Vert_{L^2}})
$$
The main result is that for $t>\tau_\mathcal{B}(u_0)$ assumption (\ref{KNass})  is fulfilled and $\Phi^t_\mathcal{B}(u_0)$ has the Koopman decomposition (\ref{finalKF}). The proof is  given in section \ref{proof4ter}.
\medskip

\noindent To make explicit the link with neighbourhoods of the sink $s_0$, let:

\begin{equation}\label{L2H1}\forall \theta_0>0\quad \Omega_\mathcal{B}(\theta_0)=\{ u_0\in L^2([0,1]);\quad  \Vert H(u_0)-\tilde{f}_0 \Vert_{L^2}<{m_0^2\over h_1(\theta_0)}\}
\end{equation}
These are neighbourhoods of the sink $s_0$ of the Burgers flow, because $H$ is continuous and the $H^1$ norm dominates the $L^2$ norm. For any Cauchy data  $u_0 \in \Omega_\mathcal{B}(\theta_0)$ the orbit  $\Phi_\mathcal{B}^t(u_0)$  has a Koopman decomposition for all $t>\theta_0$. This is proved in section \ref{proof4ter}.
\medskip

\noindent Notice that the only non-decreasing component in the decomposition cor\-responds to $(q_0,m_0)=(0,0)$. It is the sink $s_0(x)$, as expected.
%-----------------------------------------------
\subsection{on absolute convergence of formula (\ref{finalKF})}\label{ssAKB} 
Using formula (\ref{finalKF}) to compute the Koopman decomposition of observables such as  mappings of  the state variable $u(t,.)$ through analytic functions, or observables like the energy  that shows products and spatial integration,  needs convergence of (\ref{finalKF}) in the $\sup_{x\in [0,1]}$ norm. 
    
\noindent In order to have  a workable  assumption that grants absolute convergence, let:
$$
h_5(t)=\sqrt{\sum_{q=1}^\infty  {e^{-2(\mu_{q}-\mu_0)t}}(1+\sqrt{\mu_q})^2}
$$
It is a continuous decreasing function with $h_5(0)=\infty$ and $h_5(\infty)=0$. For  $u_0 \in L^2$ let $\tilde{\tau}_\mathcal{B}$ be defined by:
$$h_5(\tilde{\tau}_\mathcal{B}(u_0))={m_0^2\over \Vert H(u_0)-c_0(H(u_0))e_0)\Vert_{L^2}}
$$
The main result is that for $t>\tilde{\tau}_\mathcal{B}(u_0)$ the Koopman decomposition (\ref{finalKF}) of $\Phi^t_\mathcal{B}(u_0)$ is absolutely convergent.
\medskip

\noindent The proof is given in (\ref{proof5}).
\medskip

\noindent An interesting interpretation of this result follows  orthonormality of the fa\-mily $(e_q)$ that gives:
$$\Vert H(u_0)-\tilde{f}_0\Vert^2_{L^2}=\Vert H(u_0)-c_0(H(u_0))e_0\Vert^2_{L^2}+\Vert c_0(H(u_0))e_0-\tilde{f}_0\Vert^2_{L^2}
$$
so $$\Vert H(u_0)-c_0(H(u_0))e_0\Vert_{L^2}<\Vert H(u_0)-\tilde{f}_0\Vert_{L^2}$$
therefore if:
\begin{equation}\label{finalKB}\tilde{\Omega}_\mathcal{B}(\theta_0)=\{ u_0\in L^2([0,1]);\quad \Vert H(u_0)-\tilde{f}_0\Vert_{L^2}<     { m_0^2  \over h_5(\theta_0)}  \}
\end{equation}
the result asserted above implies that the Koopman decomposition is absolutely convergent for  $t>\theta_0$ for  Cauchy data $u_0 \in \tilde{\Omega}_\mathcal{B}(\theta_0)$.

\noindent The sets $\tilde{\Omega}_\mathcal{B}(\theta_0)$ are neighbourhoods of the sink $s_0$ in $L^2([0,1])$, by continuity of $H$ and because on $\mathcal{P}_1$ the $H^1$ norm dominates the $L^2$ norm.

%-----------------------------------------------
\subsection{Koopman eigen-observables for the Burgers flow}\label{ssKB} 
Let  $\mathcal{F}$ be the set of continuous maps from $L^2([0,1])$ to ${\mathbb{R}}$. The Koopman flow ${\mathcal{K}}^t_{\mathcal{B}}$ associated with the Burgers flow $\Phi^t_\mathcal{B}$ is the flow on  $\mathcal{F}$ given by:
\begin{equation}
    \label{fKB}
\forall \phi \in \mathcal{F}, \forall t\geq0, \forall u_0 \in L^2(]0,1[), \quad ({\mathcal{K}}^t_{\mathcal{B}}(\phi))(u_0)=\phi(\Phi^t_{\mathcal{B}}(u_0))
\end{equation}
Note that for any $t\geq0$, ${\mathcal{K}}^t_{\mathcal{B}}$ is linear on $\mathcal{F}$ and fulfills the multiplicative property:
\begin{equation}
    \label{mult}
\forall t\geq 0,\forall \phi_1 \in \mathcal{F}, \forall \phi_2 \in \mathcal{F}, \quad {\mathcal{K}}^t_{\mathcal{B}}(\phi_1 \phi_2)= {\mathcal{K}}^t_{\mathcal{B}}(\phi_1) {\mathcal{K}}^t_{\mathcal{B}}( \phi_2)
\end{equation}
Formula (\ref{iBN}) intertwins Burgers flow with the nonlinear heat flow. Formulas (\ref{vpKN2})  and (\ref{vpKN3}) give eigen-observables for the nonlinear heat flow. This gives eigen-observables for the Burgers flow as follows: 

$$\forall m\in {\mathbb{N}}, \forall \nu=(q_0,..,q_m)\in ({\mathbb{N}\backslash \{ 0\}})^{m+1},\quad  \forall u_0\in L^2([0,1])
$$
$$(\mathcal{K}^t_\mathcal{B}(\psi_\nu \circ H))(u_0)=(\psi_\nu \circ H)((\Phi^t_\mathcal{B}(u_0))=\psi_\nu((H\circ \Phi^t_\mathcal{B})(u_0))=\psi_\nu((\Phi^t_\mathcal{N}\circ H)(u_0))=
$$
$$\psi_\nu(\Phi^t_\mathcal{N}(H(u_0)))=\mathcal{K}^t_\mathcal{N}(\psi_\nu)(H(u_0))=e^{-\lambda_\nu t }(\psi_\nu \circ H)(u_0)
$$
This proves, using notations of formula (\ref{finalKF}), that for any $m\in {\mathbb{N}}$, and any $ \nu=(q_0,..,q_m)\in ({\mathbb{N}\backslash \{ 0\}})^{m+1}$:
\begin{equation}\label{vpKB1}
 e^{-\lambda_\nu t }
\end{equation}
is an eigenvalue of $\mathcal{K}^t_\mathcal{B}$ with associated eigen-observable
\begin{equation}\label{vpKB2}
\varphi_\nu=\prod_{i=1}^k{c_{q_i}\circ H\over c_{0}\circ H}
\end{equation}
To conclude this section it is worth noticing that formulas (\ref{vpKB1}) and (\ref{vpKB2}) show that the eigenvalues of the Koopman decomposition (\ref{finalKF}) are  eigenvalues of the Koopman operator $\mathcal{K}_\mathcal{B}^t$ associated with $\Phi^t_\mathcal{B}$ , and that the coefficients in the Koopman decomposition (\ref{finalKF})  are values of  eigen-functionals at the Cauchy data.

\noindent Although convergence of the Koopman decomposition  (\ref{finalKF}) is a local phenomena, eigen-functionals are defined globally on $\mathcal{F}$. 
%-------------------------------------------------
\section{proofs for the estimates}\label{proofs}
%---------------------------------------------------
\subsection{proof of formula  (\ref{estimation1})}\label{proof1}
Let $v_0 \in \mathcal{P}_+$. Let $v(t,x)$ solve equation (\ref{eqH}), so it is given by formula (\ref{devH}), and $g(t)$ is given by:
$$g(t)=\sum_0^\infty c_n(v_0)e^{-\mu_n t}\int_0^1e_n(\xi)d\xi
$$
so
$$\vert e^{\mu_0 t}g(t)-c_0(v_0)\int_0^1e_0(\xi)d\xi \vert\leq e^{-(\mu_1-\mu_0) t}\sum_1^\infty\vert c_n(v_0)\vert \vert \int_0^1e_n(\xi)d\xi \vert \leq
$$
$$e^{-(\mu_1-\mu_0) t}\sqrt{\sum_1^\infty\vert c_n(v_0)\vert^2} \sqrt{\sum_1^\infty( \int_0^1e_n(\xi)d\xi)^2}=
$$
$$e^{-(\mu_1-\mu_0) t}\Vert v_0-c_0(v_0)e_0(x)\Vert_{L^2} \Vert 1-c_0(1)e_0(x)\Vert_{L^2}
$$
%-----------------------------------------------------------
\subsection{proof of formula (\ref{estimation2})}\label{proof2}
Let $v_0 \in \mathcal{P}_+$. Let $v(t,x)$ solve equation (\ref{eqH}): it is given by formula (\ref{devH}). Three inequalities are needed for the proof:

\noindent First  is the usual Sobolev injection inequality:
\begin{equation}\label{mu4}\forall f \in H^1\quad \sup_x\vert f(x)\vert\leq C_1\Vert f\Vert_{H^1}
\end{equation}
Second is the following estimate for the $H^1$ norm of $e_n(x)$:
\begin{equation}\label{mu2}\Vert e_n\Vert_{H^1}\leq \sqrt{1+\mu_n}
\end{equation}
proved through:
$${1\over \mu_n}\int_0^1(\partial_x e_n)^2\leq {1\over \mu_n}(\int_0^1(\partial_x e_n)^2+\int_0^1Ve_n^2)\leq $$
$${1\over -\mu_n}\int_0^1(\partial_{xx}^2-V)e_n(\xi).\,e_n(\xi)\,d\xi=\Vert e_n\Vert_{L^2}=1$$

\noindent Third, because $(\partial_{xx}^2-V)^{-1}$ is bounded from $H^{-1}$ to $H^{1}$, one has:
\begin{equation}\label{mu3}\sum_0^\infty{1\over \mu_n}<\infty
\end{equation}
The proof goes as follows:
$$\sup_{x\in [0,1]}\vert v(x)-c_0(v_0)e^{-\mu_0 t}e_0(x)\vert \leq C_1\Vert v(x)-c_0(v_0)e^{-\mu_0 t}e_0(x)\Vert_{H^1}=
$$
$$C_1\Vert \sum_1^\infty e^{-\mu_nt}c_n(v_0)e_n(x)\Vert_{H^1}\leq C_1 \sum_1^\infty e^{-\mu_nt}\vert c_n(v_0)\vert \Vert e_n(x)\Vert_{H^1}\leq
$$
$$C_1 \sum_1^\infty e^{-\mu_nt}\vert c_n(v_0)\vert \sqrt{1+\mu_n}\leq C_1 \sqrt {\sum_1^\infty \vert c_n(v_0)\vert^2}\sqrt{\sum_1^\infty (1+\mu_n)e^{-2\mu_nt}}=
$$
$$e^{-\mu_0t}\Vert v_0-c_0(v_0)e_0\Vert_{L^2}h_1(t)$$
with
$$
h_1(t)=C_1\sqrt{\sum_1^\infty (1+\mu_n)e^{-2(\mu_n-\mu_0)t}}
$$
The function $h_1(t)$ is defined for all $t>0$  because of the third inequality above. It is obviously decreasing, and has limit zero at infinity.
%------------------------------------------
\subsection{proof of formula (\ref{KDN})}\label{proof3}
Let $\tilde{v}_0 \in \mathcal{P}_1$. The condition $t> \tau_\mathcal{N}(\tilde{v}_0)$ implies $\vert \tilde{k}(t)\vert<1$ so 
$${1\over 1+\tilde{k}(t)}=\sum_{m=0}^\infty(-1)^m(\tilde{k}(t))^m
$$
One expresses the product $(\tilde{k}(t))^m$ as the sum of the products of the term of rank $q_1$ in the series corresponding to the first factor, times the term of rank $q_2$ of the second factor,..., till the term of rank $q_m$ in the last factor, and gets:
$${1\over 1+\tilde{k}(t)}=\sum_{m=0}^\infty(-1)^m\sum_{(q_1,..,q_m)\in (\mathbb{N}^*)^m}e^{-\sum_{i=1}^m(\mu_{q_i}-\mu_0)t}\prod_{i=1}^m{c_{q_i}(\tilde{v}_0)\over c_{0}(\tilde{v}_0)}\prod_{i=1}^mp_{q_i}$$
Using this last expression one re-writes formula (\ref{petitdevN}) as:
$$\tilde{v}(t,x)=\sum_{q_0=0}^\infty\sum_{m=0}^\infty\sum_{(q_1,..,q_m)\in (\mathbb{N}^*)^m}(-1)^me^{-\sum_{i=0}^m(\mu_{q_i}-\mu_0)t}\prod_{i=0}^m{c_{q_i}(\tilde{v}_0)\over c_{0}(\tilde{v}_0)}\prod_{i=1}^mp_{q_i}\,\,{e_{q_0}(x)\over \int_0^1e_0(\xi)d\xi}
$$
This proves formula (\ref{KDN})
%---------------------------------------
\subsection{proof of formula (\ref{KN})}\label{proof4}
Formula (\ref{KN}) is the formula (\ref{KDN}) proved above, with a change in the order of the summation. This change is allowed if the series (\ref{KN}) is absolutely convergent, meaning that the sum of the $\sup_{x\in [0,1]}$ norm of all terms is convergent. To prove absolute convergence one needs the following estimates:
\medskip

\noindent 1- Estimating $\sup_{x\in [0,1]}\vert e_q(x)\vert$: for any $(x,y)\in [0,1]^2$,
$$\vert e_q(x)\vert \leq \vert e_q(y)\vert+ \int_y^x\vert \partial_xe_q(\xi) \vert d\xi \leq \vert e_q(y)\vert+\sqrt{ \int_0^1( \partial_xe_q(\xi) )^2 d\xi }
$$
After integrating in $y$, Cauchy-Schwarz inequality gives:
$$\vert e_q(x)\vert \leq 1+\sqrt{ \int_0^1(\partial_xe_q(\xi) )^2 d\xi }\leq 1+\sqrt{ -\int_0^1 \partial^2_{xx}e_q(\xi)\,e_q(\xi)  d\xi }=
$$
$$1+\sqrt{ \mu_q\int_0^1 e_q^2(\xi)  d\xi -\int_0^1V(\xi)e_q^2(\xi)  d\xi}\leq1+\sqrt{ \mu_q}
$$
Therefore
\begin{equation}\label{supeq}\sup_{x\in[0,1]}\vert e_q(x)\vert \leq 1+\sqrt{ \mu_q}
\end{equation}
\medskip

\noindent 2- Estimating $\sup_{x\in [0,1]}{ \vert b_\nu(x)\vert}$: the boundary conditions give:
$$\int_0^1e_q(\xi)d\xi=
-{1\over \mu_q}\int_0^1\partial^2_{xx}e_q(\xi)d\xi + {1\over \mu_q}\int_0^1V(\xi)e_q(\xi)d\xi 
=  {1\over \mu_q}\int_0^1V(\xi)e_q(\xi)d\xi 
$$
therefore, using (\ref{supeq}), one has for any $\nu=(q_0,..,q_m)\in \mathcal{B}_m$:
$$\sup_{x\in [0,1]}\vert b_\nu(x)\vert=  \sup_{x\in [0,1]}\vert e_{q_0}(x)\vert {\prod_1^m\vert \int_0^1e_{q_i}\vert \over (\int_0^1e_0)^{m+1}}    
$$
\begin{equation}\label{bnu}
\leq {{1+\sqrt{\mu_{q_0}}\over {(\int_0^1e_0(\xi)d\xi)^{m+1} }}}
\prod_1^m{\vert \int_0^1V(\xi)e_{q_i}(\xi)d\xi \vert \over \mu_{q_i}}
\end{equation}
Moreover
$$\sum_0^\infty(\int_0^1V(\xi)e_q(\xi)d\xi )^2=\Vert V\Vert^2_{L^2}<\infty
$$
\medskip

\noindent 3- Estimating $\psi_\nu(\tilde{v}_0)$:

\noindent Let $\tilde{v}_0 \in \mathcal{P}_+$. It has a square integrable derivative, and fulfills Neumann conditions on the boundary. For $q\in \mathbb{N},  q \not= 0$, one first notice that due to orthogonality of the eigenfunctions:
$$c_q(\tilde{v}_0)=\int_0^1\tilde{v}_0(\xi)e_q(\xi)d\xi=\int_0^1(\tilde{v}_0(\xi)-c_0(\tilde{v}_0)e_0(\xi))e_q(\xi)d\xi=
$$
$${1\over \mu_q}\int_0^1\partial_x(\tilde{v}_0-c_0(\tilde{v}_0)e_0)(\xi)\partial_xe_q(\xi)d\xi+{1\over \mu_q}\int_0^1(\tilde{v}_0-c_0(\tilde{v}_0)e_0)(\xi)V(\xi)e_q(\xi)d\xi
$$
so 
$$\vert c_q(\tilde{v}_0)\vert\leq {1\over \mu_q}\Vert\partial_x(\tilde{v}_0-c_0(\tilde{v}_0)e_0) \Vert_{L^2}\Vert \partial_xe_q \Vert_{L^2}+{1\over \mu_q} (\sup{V})\Vert (\tilde{v}_0-c_0(\tilde{v}_0)e_0) \Vert_{L^2}
$$
Integration by parts gives $\Vert \partial_xe_q\Vert_{L^2}\leq \sqrt{\mu_q}\,$, so:
$$\vert c_q(\tilde{v}_0)\vert\leq {1\over \sqrt{ \mu_q}}\Vert\partial_x(\tilde{v}_0-c_0(\tilde{v}_0)e_0) \Vert_{L^2}+{1\over \mu_q} (\sup{V})\Vert \tilde{v}_0-c_0(\tilde{v}_0)e_0 \Vert_{L^2}
$$
For $C_V=\max{(1,\sup{V})}$, the above estimate writes:
\begin{equation}\label{cq}\vert c_q(\tilde{v}_0)\vert \leq {C_V\over \sqrt{\mu_q}}\Vert \tilde{v}_0-c_0(\tilde{v}_0)e_0 \Vert_{H^1} 
\end{equation}
Therefore, if $\nu=(q_0,..,q_m)\in \mathcal{B}_m$:
\begin{equation}\label{psinu}
\vert \psi_\nu (\tilde{v}_0)\vert \leq C_V^{m+1}{\Vert \tilde{v}_0-c_0(\tilde{v}_0)e_0 \Vert_{H^1}^{m+1}\over  (c_0(\tilde{v}_0))^{m+1}  }(\prod_0^m\mu_{q_i})^{-{1\over 2}}\quad {\rm if}\quad q_0\not= 0
\end{equation}
and
$$\vert \psi_\nu (\tilde{v}_0)\vert \leq C_V^{m}{\Vert \tilde{v}_0-c_0(\tilde{v}_0)e_0 \Vert_{H^1}^{m}\over  (c_0(\tilde{v}_0))^{m}  }(\prod_1^m\mu_{q_i})^{-{1\over 2}}\quad {\rm if}\quad q_0= 0
$$
\medskip
\noindent 4- Final step - convergence in the $\sup_{x\in[0,1]}$ norm: 

\noindent The series below is a series of positive numbers. Therefore convergence (and sum)  do not depend on the order in indexing, by  use of Fubini's theorem. This gives, with $v_q=\vert \int_0^1V(\xi)e_{q}(\xi)d\xi \vert$:
$$\sum_\mathcal{B}e^{-\lambda_\nu t}\vert \psi_\nu(\tilde{v}_0)\vert \sup_{x\in[0,1]}\vert b_\nu(x)\vert\leq 
$$
$$\sum_{m=0}^\infty \sum_{\nu \in \mathcal{B}_m;q_0=0}e^{-\lambda_\nu t}(1+\sqrt{\mu_{0}}){C_V^{m}\Vert \tilde{v}_0-c_0(\tilde{v}_0)e_0 \Vert_{H^1}^{m}\over {(c_0(\tilde{v}_0))^{m}  {(\int_0^1e_0(\xi)d\xi)^{m+1} }}}
{\prod_{i=1}^m v_{q_i} \over\prod_{i=1}^m \mu^{3\over 2}_{q_i}} \, +
$$
$$\sum_{m=0}^\infty \sum_{\nu \in \mathcal{B}_m;q_0\not =0}e^{-\lambda_\nu t}{(1+\sqrt{\mu_{q_0}})\over \sqrt{\mu_{q_0}}}{C_V^{m+1}\Vert \tilde{v}_0-c_0(\tilde{v}_0)e_0 \Vert_{H^1}^{m+1}\over {(c_0(\tilde{v}_0))^{m+1}  {(\int_0^1e_0(\xi)d\xi)^{m+1} }}}
{\prod_{i=1}^m v_{q_i}\over \prod_{i=1}^m\mu^{3\over 2}_{q_i}}=
$$
$$\sum_{m=0}^\infty  (1+\sqrt{\mu_{0}}){C_V^{m}\Vert \tilde{v}_0-c_0(\tilde{v}_0)e_0 \Vert_{H^1}^{m}\over {(c_0(\tilde{v}_0))^{m}  {(\int_0^1e_0(\xi)d\xi)^{m+1} }}}
(\sum_{q=1}^\infty e^{-(\mu_q-\mu_0) t}{v_{q} \over \mu^{3\over 2}_{q}} )^m\, +
$$
$$\sum_{m=0}^\infty \sum_{q_0=1}^\infty e^{-(\mu_{q_0}-\mu_0) t}{1+\sqrt{\mu_{q_0}}\over \sqrt{\mu_{q_0}}}{C_V^{m+1}\Vert \tilde{v}_0-c_0(\tilde{v}_0)e_0 \Vert_{H^1}^{m+1}\over {(c_0(\tilde{v}_0))^{m+1}  {(\int_0^1e_0(\xi)d\xi)^{m+1} }}}
(\sum_{q=1}^\infty e^{-(\mu_q-\mu_0) t}{v_{q}\over \mu^{3\over 2}_{q}})^m\leq
$$
$$\sum_{m=0}^\infty  (1+\sqrt{\mu_{0}}){C_V^{m}\Vert \tilde{v}_0-c_0(\tilde{v}_0)e_0 \Vert_{H^1}^{m}\over {(c_0(\tilde{v}_0))^{m}  {(\int_0^1e_0(\xi)d\xi)^{m+1} }}}
\Vert V \Vert_{L^2}^m(h_2(t))^m\, +
$$
$$\sum_{m=0}^\infty \sum_{q_0=1}^\infty e^{-(\mu_{q_0}-\mu_0) t}{1+\sqrt{\mu_{q_0}}\over \sqrt{\mu_{q_0}}}{C_V^{m+1}\Vert \tilde{v}_0-c_0(\tilde{v}_0)e_0 \Vert_{H^1}^{m+1}\over {(c_0(\tilde{v}_0))^{m+1}  {(\int_0^1e_0(\xi)d\xi)^{m+1} }}}
\Vert V \Vert_{L^2}^m(h_2(t))^m
$$
Formula (\ref{tautilde}) gives for $t> \tilde{\tau}_\mathcal{N}(\tilde{v}_0)$:
$$\varepsilon(t)={C_Vh_2(t)\Vert V \Vert_{L^2}\Vert \tilde{v}_0-c_0(\tilde{v}_0)e_0 \Vert_{H^1}\over c_0(\tilde{v}_0)  (\int_0^1e_0(\xi)d\xi)}<1
$$
so
$$\sum_\mathcal{B}e^{-\lambda_\nu t}\vert \psi_\nu(\tilde{v}_0)\vert \sup_{x\in[0,1]}\vert b_\nu(x)\vert\leq 
{C_0\over 1-\varepsilon(t) }\sum_{q_0=0}^\infty {e^{-(\mu_{q_0}-\mu_0) t}{1+\sqrt{\mu_{q_0}}\over \sqrt{\mu_{q_0}}}  } $$
with $C_0={1\over \int_0^1e_0(\xi)d\xi}\max ({\max{(1,\sqrt{\mu_0})} ,{\max{(1,{C_V\Vert_{L^2}\Vert \tilde{v}_0-c_0(\tilde{v}_0)e_0 \Vert_{H^1}\over c_0(\tilde{v}_0)})}}  )}$
\medskip

\noindent The last series do converge for all $t>0$ due to formula (\ref{mu3}).
%-------------------------------------------
\subsection{on formula (\ref{steadyB}): the sink for Burgers equation}\label{proofsink}
Let $u_0\in L^2([0,1])$. Let $v_0=H(u_0)$ and $v(t,.)=\Phi^t_\mathcal{H}(v_0)$. Formulas (\ref{iBN}) and (\ref{lNC}) state that 
$$\Phi^t_\mathcal{B}(u_0)=C(\Phi^t_\mathcal{N}(v_0))=C({\Phi^t_\mathcal{H}(v_0)\over g(t)})
$$
Formula (\ref{estimation1}) gives:
 $$e^{\mu_0 t}g(t)=c_0(v_0)\int_0^1e_0(\xi)d\xi+\varepsilon_1(t)$$
  where 
  $$\vert \varepsilon_1(t)\vert < e^{-(\mu_1-\mu_0)t}\Vert v_0-c_0e_0(x)\Vert_{L^2} \Vert 1-c_0(1)e_0\Vert_{L^2}$$
Formula (\ref{estimation2}) gives:

$$ e^{\mu_0t}v(t,x)=c_0(v_0)e_0(x)+\varepsilon_2(t,x)$$
 where
$$ \sup_{x\in [0,1]} \vert \varepsilon_2(t,x) \vert \leq   \Vert v_0-c_0(v_0)e_0\Vert_{L^2}h_1(t)$$
and $h_1(t)$ has limit zero at infinity.

\noindent Continuity of $C$ implies:
$$\lim_{t\rightarrow \infty} \Phi^t_\mathcal{B}(u_0)(t,x)=C(\lim_{t\rightarrow \infty}{ v(t,x)\over g(t) })=C({ e_0(x)\over \int_0^1e_0(\xi)d\xi})=C(\tilde{f}_0)=s_0
$$
%--------------------------------------
\subsection{proof for formula (\ref{finalKF})}\label{proof4bis}
Let $u_0\in L^2([0,1])$ and $v_0=H(u_0)$. Let $u(t,x)$ solve (\ref{burgers}) and $v(t,x)$ solve (\ref{eqH}). It is proved in \ref{ssKBf} that under the assumption  (\ref{KNass}) one has
$$u(t,.)=-2{e^{\mu_0 t}\partial_xv\over {c_0(v_0)e_0}}\sum_{m=0}^\infty (-1)^m\big({e^{\mu_0 t}v-c_0(v_0)e_0)\over c_0(v_0)e_0)}\big)^m
$$
By formula (\ref{devH}):
$$v(t,x)=\sum_0^\infty e^{-\mu_q t}c_q(v_0)e_q(x)
$$
so
$$
u(t,x)=\sum_{q_0=0}^\infty \sum_{m=0}^\infty 2(-1)^{m+1}e^{-(\mu_{q_0}-\mu_0)t}{c_{q_0}(v_0)\over c_0(v_0)}
{\partial_xe_{q_0}(x)\over e_0(x)}\big(  \sum_{q=1}^\infty e^{-(\mu_q-\mu_0)t}{c_q(v_0)\over c_0(v_0)}{e_q(x)\over e_0(x)}\big)^m
$$
Writing the power $m$ of the series as a sum of products of the term of rank $q_1$ from the first factor,  term of rank $q_2$ from the second factor,..., till   term of rank $q_m$ from the $m$-th factor, one gets:
$$u(t,x)=$$
$$\sum_{q_0=0}^\infty \sum_{m=0}^\infty \sum_{(q_1,..,q_m)\in (\mathbb{N}\backslash \{ 0\})^m}^\infty  2(-1)^{m+1}e^{- \sum_{i=0}^m(\mu_{q_i}-\mu_0)t}
{{\prod_{i=0}^m{c_{q_i}(v_0)}\over {(c_0(v_0)})^{m+1}}}\,
{     { \partial_xe_{q_0}\prod_{i=1}^m e_{q_i}(x) }      \over {(e_0(x))^{m+1}}  }
$$
This is  formula (\ref{finalKF}).

%---------------------------------------
\subsection{on formula (\ref{L2H1})}\label{proof4ter}
Let $u_0\in L^2([0,1])$ and $v_0=H(u_0)\in \mathcal{P}_1$. Orthonormality of the sequence $(e_q)$ gives:
$$\Vert v_0-c_0(v_0)e_0\Vert_{L^2}^2=\Vert v_0-\tilde{f}_0\Vert_{L^2}^2-\Vert c_0(v_0)e_0-\tilde{f}_0\Vert_{L^2}^2
$$
so the condition (\ref{KNass}) is fulfilled if
$$\Vert v_0-\tilde{f}_0\Vert_{L^2}< {m_0 c_0(v_0)\over h_1(t)}
$$
but
\begin{equation}\label{c0m0}c_0(v_0)=\int_0^1v_0(\xi)e_0(\xi)d\xi\geq m_0\int_0^1v_0(\xi)d\xi=m_0
\end{equation}
so the condition (\ref{KNass}) is fulfilled if
$$h_1(t)< {m_0^2\over {\Vert v_0-\tilde{f}_0\Vert_{L^2}}}
$$
Because $h_1(t)$ is a decreasing continuous function with $h_1(0)=\infty$ and $h_1(\infty)=0$ this proves that for any $u_0\in L^2([0,1])$,  if
$$t>\tau_\mathcal{B}(u_0)=h_1^{-1}\big({m_0^2\over {\Vert H(u_0)-\tilde{f}_0\Vert_{L^2}}}\big)
$$
then the decomposition (\ref{finalKF}) holds.

\noindent Continuity of the map $H$ from $L^2([0,1])$ to $\mathcal{P}_1$ equipped with the $L^2$ norm shows that the map $\tau_\mathcal{B}$ is continuous from $L^2([0,1])$ to $\mathbb{R}_+$.

%---------------------------------------
\subsection{on absolute convergence of formula (\ref{finalKF})}\label{proof5}
Let $u_0\in L^2(0,1])$. This section is devoted to prove  absolute convergence in the $\sup_{[0,1]}$ norm of the following series, i.e. of the Koopman decomposition of $\Phi^t_\mathcal{B}$):
$$u(t,x)=\sum_\mathcal{B} e^{-\lambda_\nu t}\varphi_\nu(u_0)\,a_\nu(x)
$$
that means  convergence of 
\begin{equation}\label{finalKFB}\sum_\mathcal{B} e^{-\lambda_\nu t}\vert \varphi_\nu(u_0)\vert \,\sup_{x\in [0,1]} \vert a_\nu(x)\vert 
\end{equation}
where for $\nu=(q_0,..,q_m) \in \mathcal{B}  $: $\quad\lambda_\nu=\sum_0^m(\mu_{q_i}-\mu_0)$
\begin{equation}\label{finuanuB} \varphi_\nu(u_0)=\prod_0^m{ c_{q_i}(H(u_0))\over c_0(H(u_0))}\quad \quad a_\nu(x)=2(-1)^{m+1}{\partial_xe_{q_0}(x)\over e_0(x)}\prod_1^m{e_{q_i}(x)\over e_0(x)}
\end{equation}
All estimates needed are already proven, but one: an estimate with respect to $\mu_q$ of $\sup_{[0,1]}\vert \partial_xe_q(x)\vert$: because $ \partial_xe_q(x)$ fulfills Dirichlet boundary condition one has:
$$\sup_{[0,1]}\vert \partial_xe_q(x)\vert^2\leq \int_0^1(\partial_{xx}^2e_q(\xi))^2d\xi=\int_0^1(V(\xi)-\mu_q)^2e^2_q(\xi))d\xi\leq \sup_{[0,1]}\vert V(x)-\mu_q\vert^2
$$
so
\begin{equation}\label{deq}\sup_{[0,1]}\vert \partial_xe_q(x)\vert\leq C_V(1+\mu_q)
\end{equation}
1- estimate for $\sup_{[0,1]}\vert a_\nu(x)\vert$:

\noindent Formulas (\ref{supeq}) and (\ref{deq}) give for $\nu=(q_0,..,q_m)$:
 $$\sup_{x\in [0,1]}\vert a_\nu(x)\vert \leq {2C_V\over m_0}(1+{\mu_{q_0})}\prod_{i=1}^m {  1+\sqrt{\mu_{q_i}}\over m_0 }
 $$
2- estimate for $\vert \varphi_\nu(u_0)\vert$: 

\noindent Orthonormalisation of $e_q$ gives for $q\not= 0$:
$$c_q(H(u_0)) = \int_0^1H(u_0)(\xi)e_q(\xi)d\xi =
\int_0^1(H(u_0)(\xi)-c_0(H(u_0))e_0(\xi))e_q(\xi)d\xi =
$$
$$c_q(H(u_0)-c_0(H(u_0))e_0):=c'_q(u_0)
 $$
 This, with formula (\ref{c0m0}) gives
$$\vert \varphi_\nu(u_0)\vert \leq { \prod_{i=0}^m\vert c'_{q_i}\vert \over  m_0^{m+1}}\quad{\rm if}\quad q_0\not= 0
$$
$$\vert \varphi_\nu(u_0)\vert \leq {\vert c_0(H(u_0))\vert  \prod_{i=1}^m\vert c'_{q_i}\vert \over  m_0^{m+1}}\quad{\rm if}\quad q_0= 0
$$
Therefore the generic term in formula (\ref{finalKFB}) can be estimated by:
$${2C_V\over m_0^{2(m+1)}} e^{-\lambda_\nu t}(1+{\mu_{q_0}}){\prod_{i=0}^m \vert c'_{q_i}\vert  \prod_{i=1}^m(1+\sqrt{\mu_{q_i}}) }
$$
with $c'_0(u_0)=c_0(H(u_0))$.

\noindent Because it is a positive series, the sum do not depend on the summation order. Summing first in $(q_1,..,q_m)$ one gets:
$$\sum_\mathcal{B} e^{-\lambda_\nu t}\vert \varphi_\nu(u_0)\vert \,\sup_{x\in [0,1]} \vert a_\nu(x)\vert 
\leq $$
$${2C_V\over m_0^2}\sum_{q_0=0}^\infty \vert c'_{q_0}\vert (1+{\mu_{q_0}})e^{-(\mu_{q_0}-\mu_0)t}\sum_{m=0}^\infty
\big(  \sum_{q=1}^\infty  {e^{-(\mu_{q}-\mu_0)t}\over m^2_0}(1+\sqrt{\mu_q})\vert c'_q\vert \big)^m
$$
If the following assumption is fulfilled:
\begin{equation}\label{hyp}
h_3(t,u_0)= \sum_{q=1}^\infty  {e^{-(\mu_{q_0}-\mu_0)t}\over m^2_0}(1+\sqrt{\mu_q})\vert c'_q\vert <1
\end{equation}
one gets
$$\sum_\mathcal{B} e^{-\lambda_\nu t}\vert \varphi_\nu(u_0)\vert \,\sup_{x\in [0,1]} \vert a_\nu(x)\vert 
\leq {2C_V\over m_0^2(1-h_3)}\sum_{q_0=0}^\infty \vert c'_{q_0}\vert (1+{\mu_{q_0}})e^{-(\mu_{q_0}-\mu_0)t}
$$
$$
\leq {2C_V\over m_0^2(1-h_3(t,u_0))}\sqrt{\sum_{q_0=0}^\infty \vert c'_{q_0}\vert ^2} \sqrt{\sum_{q_0=0}^\infty(1+{\mu_{q_0}})^2e^{-2(\mu_{q_0}-\mu_0)t}}
$$
$$
\leq {2C_Vh_4(t)\over m_0^2(1-h_3(t,u_0))}\sqrt{\sum_{q_0=0}^\infty \vert c'_{q_0}\vert^2} 
$$
with
$$h_4(t)=\sqrt{\sum_0^\infty (1+\mu_{q})^2e^{-2(\mu_{q_0-\mu_0}})t}
$$
The function $h_4(t)$ is a continuous function for $t>0$ following (\ref{mu3}). The full  series is convergent because:
$$\sum_{q_0=0}^\infty \vert c'_{q_0}\vert ^2 =\Vert H(u_0)-c_0(H(u_0))e_0\Vert^2_{L^2}<\infty
$$
This completes the proof, and gives absolute convergence, provided assumption (\ref{hyp}) is fulfilled.

\noindent In order to fit assumption (\ref{hyp}) in a topological setting, one notices that:
$$m_0^2 h_3(t,u_0)\leq \sqrt{\sum_{q=1}^\infty  {e^{-2(\mu_{q}-\mu_0)t}}(1+\sqrt{\mu_q})^2}\sqrt{\sum_{q=1}^\infty \vert c'_q\vert^2}=
$$
$$
h_5(t)\Vert H(u_0)-c_0(H(u_0))e_0)\Vert_{L^2}
$$
with
$$
h_5(t)=\sqrt{\sum_{q=1}^\infty  {e^{-2(\mu_{q}-\mu_0)t}}(1+\sqrt{\mu_q})^2}
$$
 $h_5(t)$ is a continuous function for $t>0$ because of (\ref{mu3}). Therefore the assupmtion  (\ref{hyp}) is fulfilled for  $t>\tilde{\tau}_\mathcal{B}(u_0)$ with $\tilde{\tau}_\mathcal{B}(u_0)$ defined by:
$$h_5(\tilde{\tau}_\mathcal{B}(u_0))={m_0^2\over \Vert H(u_0)-c_0(H(u_0))e_0)\Vert_{L^2}}
$$
%---------------------------------------------
\section{Bibliography}\label{bibl}

\noindent [1]  B. O. Koopman, Hamiltonian systems and transformation in hilbert space, Proceedings of the National Academy of Sciences
298 17, 315 (1931).

\noindent [2]  M.Balabane-M.Mendez-S.Najem, On Koopman Operator for Burgers' Equation. Phys Rev. Fluids - to appear 2022.

\noindent [3] N. J. Kutz, S. L. Brunton, B. W. Brunton, and J. L. Proctor, Dynamic Mode Decomposition: Data-Driven Modeling of
Complex Systems (SIAM, 2016).

\noindent [4] J. Page and R. R. Kerswell, Koopman analysis of burgers equation, Physical Review Fluids 3, 10.1103
uids.3.071901 (2018).

\noindent [5]  I. Mezic, Analysis of 
fluid 
flows via spectral properties of the koopman operator, Annual Review of Fluid Mechanics 45,
 357 (2013).

\noindent [6] C. W. Rowley and S. T. Dawson, Model reduction for 
flow analysis and control, Annual Review of Fluid Mechanics 49,
 387 (2017). 

\noindent [7] P. J. Schmid, Dynamic mode decomposition of numerical and experimental data, Journal of Fluid Mechanics 656, 5 (2010).

%------------------------------------------

\end{document}